\documentclass[12pt]{article}

\usepackage{amssymb}
\usepackage{amsmath}
\usepackage[usenames]{color}
\usepackage{graphicx}
\usepackage[colorlinks=true,
linkcolor=webgreen,
filecolor=webbrown,
citecolor=webgreen]{hyperref}

\definecolor{webgreen}{rgb}{0,.5,0}
\definecolor{webbrown}{rgb}{.6,0,0}
\definecolor{red}{rgb}{1,0,0}

\usepackage{color}

\newtheorem{theorem}{Theorem}

\newtheorem{lemma}[theorem]{Lemma}

\newtheorem{definition}[theorem]{Definition}

\newtheorem{conjecture}[theorem]{Conjecture}
\newtheorem{question}[theorem]{Question}

\newcommand\blfootnote[1]{%
  \begingroup
  \renewcommand\thefootnote{}\footnote{#1}%
  \addtocounter{footnote}{-1}%
  \endgroup
}

\begin{document}

\title{Gaps between divisible terms in $a^2 (a^2 + 1)$}

\date{}
\author{Tsz Ho Chan}
\maketitle

\begin{abstract}
Suppose $a^2 (a^2 + 1)$ divides $b^2 (b^2 + 1)$ with $b > a$. In this paper, we improve a previous result and prove a gap principle, without any additional assumptions, namely $b \gg a (\log a)^{1/8} / (\log \log a)^{12}$. We also obtain $b \gg_\epsilon a^{15/14 - \epsilon}$ under the abc conjecture. 
\end{abstract}

\section{Introduction and Main Results} \label{s1}

We are interested in increasing sequence of positive integers $(a_n)_{n \ge 0}$ with each term dividing the next one (i.e. $a_n | a_{n+1}$). A simple example is $(2^n)_{n \ge 0}$. Another example is $(n!)_{n \ge 1}$. These are simple recursively defined sequences. It is more interesting and challenging to require each term of the sequence to have a special form. For example, $(3^n)_{n \ge 0}$ has all terms odd. For another special example, let us consider the Fibonacci numbers
\[
F_1 = 1, F_2 = 1, F_{n+1} = F_{n} + F_{n-1} \text{ for } n \ge 2.
\]
By the well-known fact that $F_m | F_n$ if and only if $m | n$, the sequence $(F_{2^n})_{n \ge 0}$ has all terms of Fibonacci-type and each term dividing the next one. One may restrict the sequence to numbers of the form $n^2 + 1$ or other polynomials, and we are interested in the growth of such sequences. In this paper, we shall focus on numbers of the form $n^2 (n^2 + 1)$ and study the following
\begin{question} \label{ques1}
Suppose $a^2 (a^2 + 1)$ divides $b^2 (b^2 + 1)$. Must it be true that there is some gap between $a$ and $b$? More preciously, is it true that $b > a^{1 + \epsilon}$ for some small $\epsilon > 0$?
\end{question}
In [\ref{C1}], the author studied the above question with some additional restrictions on $a$ and $b$. In this paper, we remove all these restrictions and prove
\begin{theorem} \label{thm1}
Let $a$ and $b$ be positive integers with $a < b$. Suppose $a^2 (a^2 + 1)$ divides $b^2 (b^2 + 1)$. Then
\[
b \gg \frac{a (\log a)^{1/8}}{(\log \max(e, \log a))^{12}}.
\]
\end{theorem}

Recently, Stephen Choi, Peter Lam and the author [\ref{CCL}] defined gap principle of order $n$ for polynomials with integer coefficients as follows:
\begin{definition}
Let $n$ be a positive integer and $f(x)$ be a polynomial with integral coefficients. Consider the set of all positive integers $a_0 < a_1 < a_2 < \ldots < a_n$ such that $f(a_i)$ divides $f(a_{i+1})$ for $0 \le i \le n-1$. We say that $f(x)$ satisfies the gap principle of order $n$ if $\lim a_n / a_0 = \infty$ as $a_0 \rightarrow \infty$.
\end{definition}
Hence, Theorem \ref{thm1} implies that the polynomial $f(x) = x^2 (x^2 + 1)$ satisfies the gap principle of order $1$.

\bigskip

Assuming the abc conjecture, we can obtain a better gap result than that in [\ref{C1}], without any extra assumptions.
\begin{theorem} \label{thm2}
Let $a$ and $b$ be positive integers with $a < b$. Suppose $a^2 (a^2 + 1)$ divides $b^2 (b^2 + 1)$. Then, under the abc conjecture with any small $\epsilon > 0$,
\[
b \gg_\epsilon a^{15 / 14 - \epsilon}.
\]
\end{theorem}
Thus, this answers question \ref{ques1} in the affirmative under the abc conjecture. We hope this article would inspire readers to study questions of similar nature.

\bigskip

{\bf Some Notations} The symbol $a | b$ means that $a$ divides $b$. The notations $f(x) \ll g(x)$ and $g(x) \gg f(x)$ are all equivalent to $|f(x)| \leq C g(x)$ for some constant $C > 0$. Finally $f(x) = O_\lambda(g(x))$, $f(x) \ll_\lambda g(x)$ or $g(x) \gg_\lambda f(x)$ mean that the implicit constant $C$ may depend on $\lambda$.

\section{Proof of Theorem \ref{thm1}} \label{s2}

Since $a^2 (a^2 + 1)$ divides $b^2 (b^2 + 1)$, say
\[
t a^2 (a^2 + 1) = b^2 (b^2 + 1)
\]
for some integer $t > 1$. We may assume $t \le \log a$, for otherwise the theorem is true automatically. Let $D$ be the greatest common divisor of $a$ and $b$. Suppose $a = D x$ and $b = D y$ with $(x, y) = 1$, and let $T = (D^2 x^2 + 1, D^2 y^2 + 1)$. Then
\begin{equation} \label{eq0}
t x^2 \frac{D^2 x^2 + 1}{T} = y^2 \frac{D^2 y^2 + 1}{T}.
\end{equation}
Since $(x, y) = 1$, $x^2$ must divide $(D^2 y^2 + 1) / T$. Say $(D^2 y^2 + 1) / T = m x^2$ for some integer $m$. Then $t (D^2 x^2 + 1) / T = m y^2$, and we have
\begin{equation} \label{eq1}
t (D^2 x^2 + 1) = m T y^2
\end{equation}
and
\begin{equation} \label{eq2}
D^2 y^2 + 1 = m T x^2.
\end{equation}
Multiplying (\ref{eq1}) by $D^2$, (\ref{eq2}) by $m T$ and combining, we have
\begin{equation} \label{eq3}
[(m T)^2 - t D^4] x^2 = t D^2 + m T.
\end{equation}
Similarly, multiplying (\ref{eq1}) by $m T$, (\ref{eq2}) by $t D^2$ and combining, we have
\begin{equation} \label{eq4}
[(m T)^2 - t D^4] y^2 = t D^2 + t m T.
\end{equation}
Subtracting (\ref{eq3}) from (\ref{eq4}), we get
\begin{equation} \label{eq5}
s (y^2 - x^2) = (t - 1) m T
\end{equation}
with
\begin{equation} \label{eq6}
s = (m T)^2 - t D^4.
\end{equation}
From (\ref{eq2}), we have $(m T, D) = 1$. Hence,
\[
(s, m T) = ( (m T)^2 - t D^4, m T ) = (t D^4, m T) = (t, m T).
\]
Therefore, by combining this with (\ref{eq5}), we have $s | (t - 1) t$ and (\ref{eq6}) gives us an hyperelliptic curve
\[
Y^2 = t X^4 + s
\]
by setting $Y = m T$ and $X = D$. Let $\lambda := (t - 1) t$. By Theorem 1 of Voutier [\ref{V}] on the study of integral solutions to hyperelliptic curves using transcendental number theory, we have
\[
\max(X, Y) < e^{C_1 \lambda (\max(1, \log \lambda))^{96}}
\]
for some constant $C_1 \ge 1$. Suppose $\lambda \le \frac{1}{C_1} \frac{\log D}{(\log \max(e, \log D))^{96}}$. Then
\[
D = X < e^{C_1 \frac{1}{C_1} \frac{\log D}{(\log \max(e, \log D))^{96}} (\log \max(e, \log D))^{96}} = D
\]
which is a contradiction. Therefore, $(t - 1) t> \frac{1}{C_1} \frac{\log D}{(\log \max(e, \log D))^{96}}$ and
\begin{equation} \label{ineq0}
t \gg \frac{\sqrt{\log D}}{(\log \max(e, \log D))^{48}}.
\end{equation}
From (\ref{eq1}), we have $2 t D^2 x^2 \ge m T y^2$ and $t D^2 \gg m T$. This together with (\ref{eq3}) gives $t D^2 \gg x^2$. Hence, as $t \le \log a$ and $a = D x$, we have $\log D \gg \log a$. Therefore, (\ref{ineq0}) gives
\[
t \gg \frac{\sqrt{\log a}}{(\log \max(e, \log a))^{48}}
\]
and we have Theorem \ref{thm1} as $b^4 / a^4 \gg t$.

\section{Proof of Theorem \ref{thm2}} \label{s3}

Firstly, for any integer $n$, let $R(n) := \prod_{p | n} p$ be the radical or kernel of an integer $n$. The abc-conjecture is
\begin{conjecture}
For every $\epsilon > 0$, there exists a constant $C_\epsilon$ such that for all triples $(a,b,c)$ of coprime positive integers, with $a + b = c$, then
\[
c < C_\epsilon R(a b c)^{1 + \epsilon}.
\]
\end{conjecture}
Secondly, let us state a lemma which follows easily from the unique prime factorization of numbers.
\begin{lemma} \label{lem1}
Suppose $a | A^2$ and $a = a_1 a_2^2$ with $a_1$ squarefree. Then $a_1 a_2 | A$.
\end{lemma}

\bigskip

Basically, we follow the proof of Theorem \ref{thm1}. Using the same notation as in section \ref{s2},
\[
a = D x, \; b = D y, \text{ and } s + t D^4 = (m T)^2
\]
with $(m T, D) = 1$. Let
\[
(t, (m T)^2) = d \text{ and } d = d_1 d_2^2 \text{ with } d_1 \text{ squarefree}.
\]
By Lemma \ref{lem1}, we have $m T = d_1 d_2 S$. Then
\begin{equation} \label{abc1}
\frac{s}{d} + \frac{t}{d} D^4 = d_1 S^2.
\end{equation}
Here, the three terms are integers and pairwise relatively prime. We can apply the abc conjecture and obtain
\begin{equation} \label{abc2}
\frac{t}{d} D^4 < C_\epsilon \Bigl(\frac{(t - 1) t}{d} D d_1 S \Bigr)^{1 + \epsilon}
\end{equation}
since $s | (t-1) t$. Suppose $t \le 10 D^4$. As $d_1 \le d$, (\ref{abc1}) gives
\[
(d_1 S)^2 \ll t D^4 \text{ or } d_1 S \ll t^{1/2} D^2.
\]
Putting this into (\ref{abc2}), we have
\[
\frac{t}{d} D^4 \ll_\epsilon \Bigl(\frac{t^{5/2}}{d} D^3 \Bigr)^{1 + \epsilon}.
\]
Hence,
\begin{equation} \label{abc3}
D \ll_\epsilon t^{3/2 + 8 \epsilon}.
\end{equation}
As $s + t D^4 = (m T)^2$, $s | (t-1) t$ and $t \le D^4$, we have $m T \ll t^{1/2} D^2$. This together with (\ref{eq3}) gives
\[
x^2 \ll \frac{t D^2}{s} \text{ or } x \ll t^{1/2} D.
\]
Therefore, by (\ref{abc3}),
\[
a = D x \ll t^{1/2} D^2 \ll_\epsilon t^{7/2 + 16 \epsilon} \text{ or } t \gg_\epsilon a^{2/7 - 2 \epsilon}
\]
which gives Theorem \ref{thm2} as $b^4 / a^4 \gg t$.

\bigskip

We are left with the case $t > 10 D^4$. Suppose $y / x = b / a \le a^{1/5} = (D x)^{1/5}$. From (\ref{eq0}), we have $t^{1/4} \le \sqrt[4]{2} y / x$. Hence,
\begin{equation} \label{ineq1}
10 D^4 < t \le 2 (D x)^{4/5}.
\end{equation}
This gives
\begin{equation} \label{ineq2}
x > \frac{10}{2^{5/4}} D^4.
\end{equation}
On the other hand, (\ref{eq2}) and (\ref{eq0}) give
\[
m T \le 2 D^2 y^2 / x^2 \le 2 \sqrt{2} t^{1/2} D^2.
\]
Putting this into (\ref{eq3}), we have
\begin{equation} \label{ineq3}
x^2 \le (1 + 2 \sqrt{2}) t D^2 \text{ or } \frac{1}{1 + 2 \sqrt{2}} \Bigl(\frac{x}{D}\Bigr)^2 \le t.
\end{equation}
Combining (\ref{ineq1}) and (\ref{ineq3}), we obtain
\[
\frac{1}{1 + 2 \sqrt{2}} \Bigl(\frac{x}{D}\Bigr)^2 \le 2 (D x)^{4/5} \text{ or } x \le [2(1 + \sqrt{2})]^{5/6} D^{7/3}
\]
which contradicts (\ref{ineq2}). Therefore, we must have $b / a > a^{1/5}$ which also gives Theorem \ref{thm2}.


\noindent
Department of Mathematical Sciences \\
University of Memphis \\
Memphis, TN 38152 \\
U.S.A. \\
thchan6174@gmail.com

\end{document}